# Comparing the Performance of Optimal PID and Optimal Fractional-Order PID Controllers Applied to the Nonlinear Boost Converter

F. Merrikh-Bayat* A. Jamshidi**

*Department of Electrical and Computer Engineering, University of Zanjan, Zanjan
IRAN (e-mail: f.bayat@znu.ac.ir)
**Department of Electrical and Computer Engineering, University of Zanjan, Zanjan, IRAN (e-mail: ar_1364@yahoo.com)

**Abstract:** This paper proposes the application of fractional-order PID (FOPID) controller for output voltage control of boost converters. For this purpose, parameters of the FOPID controller are calculated such that the Integral Absolute Error (IAE) of the variations of the output voltage is minimized. Since the search space is very large in dealing with such an optimization problem, the Artificial Bee Colony (ABC) algorithm is used for optimal tuning the parameters of the FOPID controller. Simulations, which are performed by using the complete non-minimum phase model of the boost converter, confirm the fact that the proposed optimal FOPID controller can improve the transient response of the feedback system by using a considerably smaller control effort (i.e., less on-off switches) compared to the optimal PID controller. Moreover, it is shown that the proposed FOPID controller enhances the robustness of the boost converter to variations in the input voltage.

## 1. INTRODUCTION

By spreading usage of electronic devices and paying more attention to renewable sources, DC-DC converters become a more interesting field of study for more researchers. In many applications it is highly demanded to design a DC-DC converter with the ability of tracking the reference voltage with appropriate transient and steady-state responses at the presence of input voltage variations. However, since these converters are inherently time-varying and non-linear, fulfillment of the abovementioned requirements may not be possible without using complicated control techniques (Almer et al., 2007; Beccuti et al., 2005; Beccuti et al., 2006; Beccuti et al., 2007; Calderon et al., 2005; Guo et al., 2003; Perry et al. 2007; Rafiei et al., 2003; Tsang and Chan, 2005).

So far, many different methods have been proposed for the control of DC-DC converters. For example, design of cascade controller for DC-DC buck converter (Tsang and Chan, 2005), model predictive control of the DC-DC boost converter (Beccuti et al., 2006), Hybrid control of DC-DC converters (Almer et al., 2007; Beccuti et al., 2007), optimal control of the boost DC-DC converter (Beccuti et al., 2005), PI-like fuzzy logic control of DC-DC converter (Perry et al. 2007), robust control of DC-DC PWM converters (Rafiei et al., 2003), and digital PI control of DC-DC converters (Guo et al., 2003) can be found in the literature. Some of these methods can be easily designed and realized but are less effective, while some others are very effective but very difficult to realize. For example, hybrid control techniques are often (theoretically) very effective for DC-DC converters, but it is extremely difficult to realize such controllers since the computational effort needed to calculate the control signal is very high in this case and needs special hardware, which is not always within reach. Moreover, such methods need the exact model of the converter, which is not always available. Analog controllers are advantageous in the way that can generate real-time control signals. But, these controllers are mainly limited to PID-type compensators, which have a limited performance especially in dealing with non-minimum phase systems such as boost converters.

The main aim of this paper is to propose the application of a new generation of PID controllers, which is called the fractional-order PID (FOPID), for optimal output-voltage control of boost converters. It should be noted that since FOPIDs have five parameters to tune (i.e., two more than the conventional PID controllers), they can be applied to more complicated control problems. Moreover, the performance of control is often much better when these controllers are applied. It is also important to note that any improvement in the quality of the output voltage of a DC-DC converter is of high importance in practice.

According to the non-minimum phase behavior of boost DC-DC converters, and consequently, the more difficulties one may face in dealing with such converters (compared to, e.g., the buck converters), all studies of this paper are mainly focused on this type of converters. Clearly, the basic ideas of this work can easily be extended to other types of converters. To save the generality of discussions, the complete nonlinear model of the boost DC-DC converter is used in all simulations.

The rest of this paper is organized as follows. Section 2 represents the fundamentals of fractional-order control systems. The so-called fractional-order PID controller is also reviewed in this section, and a commonly-used method for the simulation of these controllers is briefly discussed. In Section 3 we will use the Artificial Bee Colony (ABC) algorithm for optimal tuning the parameters of a FOPID controller, which will be applied to a boost converter. Some notes in relation to optimal tuning the FOPID controllers are also discussed in this section. Section 4 represents a comparison between the performance of optimal PID and optimal FOPID controllers when applied to a boost converter. Finally, Section 5 concludes the paper.

## 2. REVIEW OF FRACTIONAL-ORDER PID CONTROL

Due to the simplicity in design and good performance (i.e., low overshoot and small settling time), PID controllers are of the most widely-used controllers in industry. Nevertheless, because of the limitations of this type of controller, there had always been a continuous attempt to improve the performance and robustness of PID controllers (Johnson, 2005; Wang et al., 1999; Jin et al., 1998; Shahruz and Schwartz, 1997). In recent years, according to the advances in the field of fractional calculus, there had been a great interest to develop a new generation of PID controllers, which is commonly known as the fractional-order PID (FOPID) or $PI^\lambda D^\mu$ controller. The transfer function of a FOPID controller, which was initially proposed by Podlubny in 1999, is given by

$$G_c(s) = K_p \left( 1 + \frac{1}{T_i s^\lambda} + T_d s^\mu \right) \quad (1)$$

where $K_p, T_i, T_d \in \mathbb{R}$ and $\lambda, \mu \in \mathbb{R}^+$ are the tuning parameters and the controller design problem is to determine the suitable value of these unknown parameters such that a predetermined set of control objectives is met (Podlubny, 1999a). Note that in (1) the fractional powers of the Laplace variable, $s$, are commonly interpreted in the time domain using either the Grunwald-Letnikov or the Riemann-Liouville or the Caputo definition (Podlubny, 1999b).

It should be noted that any classical PID-controller is a particular case of the FOPID controller given in (1). For example, assuming $\lambda = \mu = 1$ the FOPID controller given in (1) is reduced to the conventional PID controller. It concludes that, in general, the FOPID controller is more flexible than the PID controller and provides us with an opportunity for better adjustment of the dynamical properties of the feedback system under consideration, of course at the cost of using a more complicated setup. Various successful applications of FOPID controllers can be found in the literature (see, for example, (Ma and Hori, 2007; Monje et al., 2007; Oustaloup et al., 1995) and the references therein for more information on this subject).

As it can be observed, the FOPID controller given in (1) has five parameters to tune, i.e. two more than the conventional PID controllers, and it is the main reason for the superiority of FOPIDs to PIDs. The frequently used strategy for tuning FOPID controllers is based on the minimization of a suitably chosen cost function, commonly in the time domain. For example, optimal tuning rules of the FOPID controller, when the process is modeled with a first-order plus time delay (FOPTD) transfer function and minimization of either the ISE or the ISTE performance index is aimed are resented in (Merrikh-Bayat, 2011). Since (optimal) tuning of FOPID controllers needs searching the very large five-dimensional space, nature-inspired population-based optimization algorithms are commonly used for this purpose. For example, application of the Particle Swarm Optimization (PSO) for minimizing the Integral Time Absolute Error (ITAE) (Maiti et al., 2008), PSO to minimize a weighted combination of ITAE and control input (Cao et al., 2006), hybrid of Electromagnetism-like (EM) and Genetic Algorithm (GA) to minimize the ISE (Chang and Lee, 2007), and GA to minimize a combination of ITAE and control input (Cao et al., 2005) can be found in the literature.

Although the FOPID controllers can effectively improve the performance of control, this improvement is achieved at the cost of using a more complex system needed to realize fractional-order integrators and differentiators. So far, few methods have been developed by researchers for the simulation and realization of FOPID controllers (Charef et al., 1992; Charef, 2006; Valerio and Costa, 2005; Vinagre et al., 2000). Among others, Oustaloup recursive approximation (ORA) (Oustaloup, 1995; Oustaloup et al., 2000) is the most popular one used for this purpose. In this method, the fractional-order differentiator $s^v$ ($v > 0$) is approximated as

$$s^v \approx k \prod_{n=1}^{N} \frac{1 + s/\omega_{z,n}}{1 + s/\omega_{p,n}} \quad (2)$$

where $N, \omega_{z,n}, \omega_{p,n}$ and $k$ are the unknown parameters to be determined. ORA provides us with a set of recursive formulas for determining the value of these parameters such that the approximation be valid in the desired frequency range $[\omega_l, \omega_h]$. For this purpose, the number of poles and zeros, $N$, is chosen beforehand, and then the frequency of poles and zeros is calculated by using the following recursive equations:

$$\omega_{z,1} = \omega_l \sqrt{\eta} \quad (3)$$

$$\omega_{p,n} = \omega_{z,n} \alpha \quad n = 1, \ldots, N \quad (4)$$

$$\omega_{z,n+1} = \omega_{p,n} \eta \quad n = 1, \ldots, N-1 \quad (5)$$

where

$$\alpha = \left( \omega_h / \omega_l \right)^{v/N} \quad (6)$$

$$\eta = \left( \omega_h / \omega_l \right)^{(1-v)/N} \quad (7)$$

After calculating the frequency of poles and zeros, the gain $k$ in (2) is adjusted so that both sides of this equation have unit gains at 1 rad/s. In this method, the good performance of the approximation strongly depends on the value of $N$. More precisely, low values result in simpler approximations but also cause the appearance of a ripple in both gain and phase behaviors. Such a ripple can be eliminated by increasing the value of $N$, clearly at the cost of a more computational effort or a more complex setup (see (Merrikh-Bayat, 2012) for some hints on the selection of $N$, $\omega_l$ and $\omega_h$). The case $v<0$ can be dealt with by inverting (2). Since the approximation becomes unsatisfactory for $|v|>1$, it is common to split the fractional powers of $s$ as

$$s^v = s^n s^\delta \quad (8)$$

where $n \in \mathbb{Z}$ and $\delta \in [0,1]$. In this manner only the latter term in (8) has to be approximated. In all simulations of this paper the ORA is used to approximate the fractional-order differentiator and integrator of the FOPID controller.

## 3. OPTIMAL FOPID CONTROLLER DESIGN FOR BOOST CONVERTER

Optimal tuning the parameters of a FOPID controller, which is applied to a boost converter, first needs the definition of optimization goal(s). In this paper, the aim of optimization is to find the optimal values of $K_p, T_i, T_d, \lambda$, and $\mu$ such that the modified IAE index defined as follows

$$J_{IAE} = \int_0^T |v_{out} - v_{ref}| dt \quad (9)$$

is minimized, where $v_{out}$ and $v_{ref}$ are the output and the reference voltage of the boost converter, respectively (see Fig. 2). It is a well-known fact that using the IAE index for controller design leads to a feedback system with a small overshoot and a small settling time, which are highly desired in many control applications. But, the drawback of using this performance index is that, in general, it cannot be calculated analytically and one needs to apply numerical techniques to evaluate it. Note that because of the switching action used in all DC-DC converters, the output voltage of such converters (including the boost DC-DC converter) will not asymptotically tend to the reference voltage. Hence, the standard IAE index (or any other integral performance index) will not be convergent in this case and one should instead use a performance index with a finite period of integration (e.g., such as the one defined in (9)). Moreover, in order to arrive at a feedback system with a good transient and steady-state response, the value of $T$ in (9) should be properly chosen, i.e., it should be sufficiently larger than the time constant of the closed-loop system under consideration.

In order to calculate the parameters of the optimal FOPID controller under consideration, first we realize it in SIMULINK as shown in Fig. 2 (Table 1 contains the value of the elements used in this converter). In this figure, $G_c(s)$ is the transfer function of the FOPID controller, which is realized by using the ORA method described in Section 2.

After realizing the feedback system shown in Fig. 2 in SIMULINK, we simulate it repeatedly by using the FOPID controller whose parameters are obtained from the Artificial Bee Colony (ABC) algorithm (Karaboga and Basturk, 2007) (in fact, a m-file code generates the parameters of the FOPID controller by using the ABC algorithm and then sends them to the SIMULINK. After running the SIMULINK model, the value of the objective function is returned back to the m-file code). After each simulation, the transfer function of the best controller can be determined by comparing the value calculated for $J_{IAE}$ in the last simulation with those calculated in previous simulations. The ABC algorithm runs until the stop criterion is reached. Note that in order to find the parameters of the optimal FOPID controller any other nature-inspired optimization algorithm can also be applied.

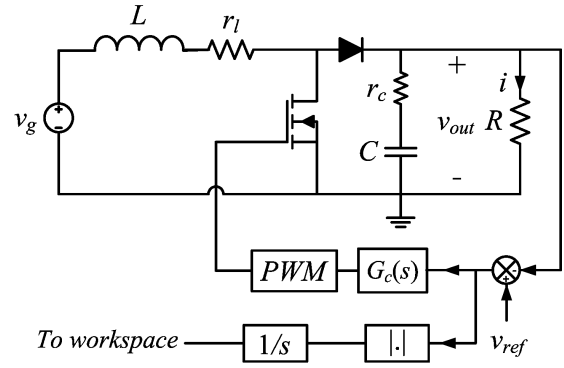

Fig. 2. SIMULINK model of the boost converter

**Table 1. Parameters of the boost converter shown in Fig. 2**

| Parameter | Value | Unit |
|---|---|---|
| Load resistance, $R$ | 25 | Ω |
| Filter inductance, $L$ | 250 | μH |
| $r_l$ | 0.075 | Ω |
| Filter capacitance, $C$ | 1056 | μF |
| $r_c$ | 0.0375 | Ω |

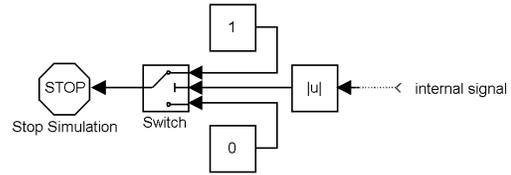

Fig. 3. The breaking mechanism used to avoid infeasible solutions (i.e., those lead to unstable feedback systems) during the numerical simulation

There is another point that should be noted during the numerical optimization and that is since the values suggested for $K_p, T_i, T_d, \lambda$, and $\mu$ by the ABC algorithm may result in an unstable feedback system, we should provide a breaking mechanism in the SIMULINK model to stop the simulation when the signals are going to blow up. Clearly, the optimization procedure will halt without applying such a breaking mechanism. One possible approach to stop the simulation when the feedback system shown in Fig. 2 is unstable is to apply the breaking system shown in Fig. 3 to all internal signals of the SIMULINK model. In this figure the switching threshold should be chosen equal to a big number (e.g., $10^{20}$).

## 4. COMPARING THE PERFORMANCE OF INTEGER-ORDER AND FRACTIONAL-ORDER PID CONTROLLERS APPLIED TO THE BOOST CONVERTER

This section presents a comparison between the results obtained by applying optimal PID and FOPID controllers to the boost converter shown in Fig. 2. As mentioned before, parameters of the optimal PID and FOPID controllers used in the following simulations are calculated using the ABC. It is also assumed that the boost converter under consideration converts the 5 V input voltage to 12 V at the output (i.e., $v_g = 5$ V and $v_{ref} = 12$ V in Fig. 2). All of the following simulations are performed assuming $T = 0.05$ sec (see (9)). Moreover, during the numerical optimization it is assumed that the ABC algorithm stops after 20 iterations. The number of artificial bees is also considered equal to 10.

Since all nature-inspired population-based optimization algorithms are of stochastic nature, a given optimization problem should be solved several times in order to make sure about finding the global best solution with a high probability. Tables 2 and 3 show the results of applying the ABC algorithm to the problems of designing optimal PID and FOPID controller, respectively. Each table contains the results of four optimizations, where the best solutions are shown by boldface characters. As it can be observed, the best value of $J_{IAE}$, when optimal PID and optimal FOPID controllers are applied, is equal to 0.0175 and 0.0168, respectively. It is also observed that the optimal FOPID controller can effectively decrease the settling time of the boost converter compared to the optimal PID controller.

Figures 4 and 5 show the output voltage of the boost converter when optimal PID and optimal FOPID controllers are applied, respectively. Although the IAE performance index defined in (9) mainly takes into account the transient response of the boost converter, but as it can observed in Figs. 4 and 5 it also improves the steady-state response of the system (notice to the values of output signals at the locations referred to by borders on these figures).

In order to make a better comparison between the performance of optimal PID and FOPID controllers, the corresponding steady-state and transient responses are shown in Figs. 6 and 7, respectively with more details. According to Fig. 6 the output voltage is more stable and better regulated when the FOPID controller is applied. Figure 7 clearly shows that the optimal FOPID controller leads to a feedback system with a considerable faster response. Interesting observation is that all of the abovementioned improvements are achieved by using even a smaller control effort. More precisely, counting the number of on-off switches in the time period $[0, 0.05]$ reveals the fact that PID and FOPID controllers apply 383 and 192 on-off switches, respectively. In other words, the FOPID controller uses 49.86% less switching actions compared to the PID controller, which can highly reduce the transient disturbances and the losses due to switching.

Figure 8 shows the effect of adding the 10% step disturbance to the input voltage (which occurs at $t = 0.15$ sec) on the output voltage when the optimal PID and the optimal FOPID controllers are applied. As it can be observed, the FOPID controller rejects the effect of this disturbance much better than the PID controller.

**Table 2. Results of applying ABC algorithm to the problem of designing optimal PID controller**

| Simulation number | 1 | 2 | 3 | 4 |
|---|---|---|---|---|
| $K_p$ | **3.6058** | 1.2174 | 1.7985 | 3.8759 |
| $T_d$ | **0.0147** | 0.0438 | 0.0060 | 0.0608 |
| $T_i$ | **2.7894e-4** | 8.1197e-5 | 5.9339e-4 | 4.7149e-5 |
| $J_{IAE}$ | **0.0175** | 0.0179 | 0.0178 | 0.0176 |
| Overshoot (%) | **0.00** | 0.00 | 0.00 | 0.00 |
| Settling time (sec) | **0.0049** | 0.0114 | 0.0131 | 0.0076 |

**Table 3. Results of applying ABC algorithm to the problem of designing optimal FOPID controller**

| Simulation number | 1 | 2 | 3 | 4 |
|---|---|---|---|---|
| $K_p$ | 1.2864 | **1.7274** | 3.3400 | 2.1678 |
| $T_d$ | 0.0449 | **0.0345** | 0.0563 | 0.0836 |
| $T_i$ | 3.7493e-4 | **9.3187e-4** | 8.1715e-4 | 2.9188e-4 |
| $\mu$ | 0.9377 | **0.8914** | 0.7255 | 0.8635 |
| $\lambda$ | 0.7493 | **0.7157** | 0.8977 | 0.8181 |
| $J_{IAE}$ | 0.0171 | **0.0168** | 0.0169 | 0.0175 |
| Overshoot (%) | 1.7240 | **0.000** | 0.000 | 0.000 |
| Settling time (sec) | 0.0066 | **0.0019** | 0.0017 | 0.0098 |

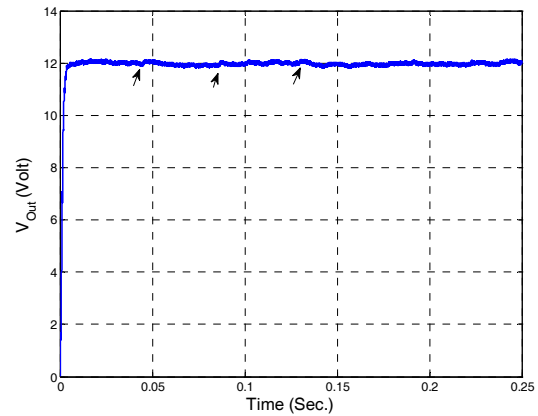

Fig. 4. Output voltage of the PID-controlled boost converter

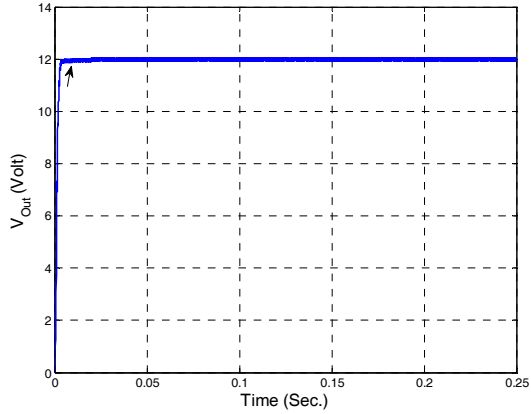

Fig. 5. Output voltage of the FOPID-controlled boost converter

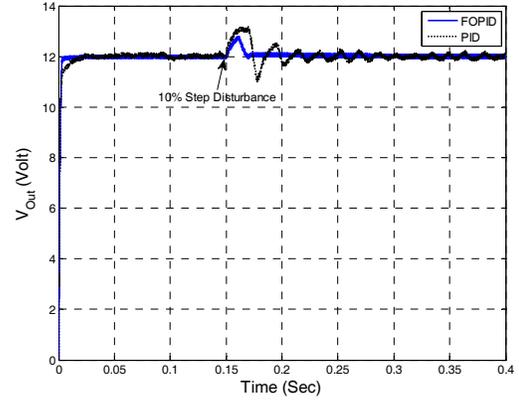

Fig. 8 The effect of adding a 10% step disturbance to the input voltage on the output voltage

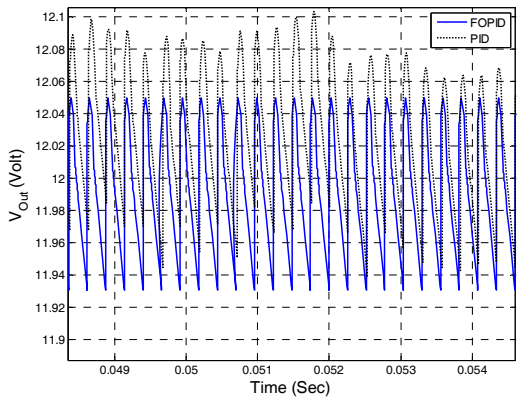

Fig. 6. Steady-state value of the output voltages shown in Figs. 4 and 5 with more details

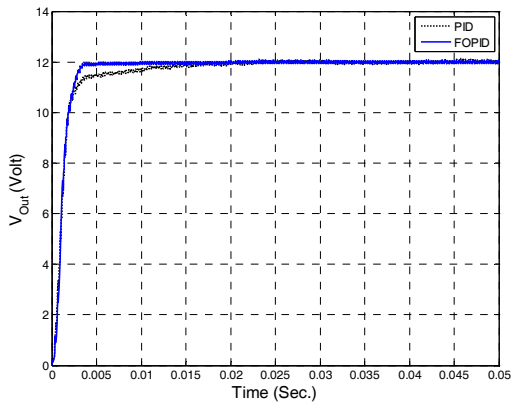

Fig. 7. Start up responses of the output voltages shown in Figs. 4 and 5 with more details

## 5. CONCLUSION

In this paper, application of the fractional-order PID (FOPID) controller for output voltage control of the boost DC-DC converter is studied and a method for optimal tuning the parameters of this controller, which is based on the artificial bee colony algorithm, is proposed.

Simulation results show that the proposed FOPID controller can improve the startup response of the boost converter by using less on-off switching actions compared to the optimal PID controller. This result is of high importance in practice since reducing the number of on-off switches can effectively decrease the transient disturbances and losses due to switching. Simulations also prove that the proposed FOPID controller can effectively improve the rejection of possible disturbances, which may occur in the input voltage. It is also observed that the output voltage is better regulated when the FOPID controller is applied.

One main advantage of the proposed method for optimal tuning the FOPID controllers (for boost DC-DC converters) is that it can work with any model of the boost converter. In fact, unlike many other methods that strictly depend on the mathematical model of the converter, the proposed method works without the need to such mathematical models.